\documentclass[11pt]{amsart}
\usepackage{latexsym}
\usepackage{amssymb,amsmath,amsthm,amscd, 
amssymb,enumerate, verbatim, stmaryrd}
\def\co{\colon\thinspace}

\newcommand{\curv}{\mbox{curv}}
\newcommand{\inj}{\mbox{inj}}

\newcommand{\diam}{\mbox{diam}}
\newcommand{\Ric}{\mbox{Ric}}

\newtheorem{thm}{Theorem}[section]

\newcommand{\Vol}{\mbox{vol}}
\newcommand{\seccurv}{\mbox{sec}}

\newtheorem{Example}[thm]{Example}

\newtheorem{remark}[thm]{Remark}

\newtheorem{Fact}[thm]{Fact}

\newtheorem{Main Lemma}[thm]{Main Lemma}

\newtheorem{Convention}[thm]{Convention}

\begin{document}\abovedisplayskip=6pt plus3pt minus3pt 
\belowdisplayskip=6pt plus3pt minus3pt
\title{\bf
Degenerations of Riemannian manifolds\rm}
\author{Igor Belegradek}
\date{}
\thanks{\it 2000 Mathematics Subject classification.\rm\ Primary
53C20.}\rm

\address{Igor Belegradek\\School of Mathematics\\ Georgia Institute of
Technology\\ Atlanta, GA 30332-0160}\email{ib@math.gatech.edu}

\maketitle

\section{Motivation and basic definitions}

This is an expositiry article on collapsing theory 
written for the Modern Encyclopedia of Mathematical Physics (MEMPhys).
We focus on describing the geometric and topological structure of
collapsed/non-collapsed regions in Riemannian manifold
under various curvature assumptions. 
Numerous applications of collapsing theory to
Riemannian geometry are {\it not} discussed in this survey, 
due to page limits dictated by the encyclopedia format.
More information on collapsing can be found in
the ICM articles by Perelman~\cite{perelman-icm},
Colding~\cite{col-icm},
Petrunin~\cite{Petrunin-icm}, Rong~\cite{Rong-icm}, 
in Cheeger's book on Cheeger-Colding theory~\cite{Chi-Ric-book}, 
and in the comprehensive survey of Fukaya~\cite{Fuk2006}.
The article ends with
an appendix on Gromov-Hausdorff distance, also written for  
MEMPhys.

A fundamental problem in Riemannian geometry is to analyze
how a family of Riemannian manifolds can degenerate.
For example, in the case of Einstein manifolds it is natural 
to study how the Einstein equation develops a singularity, or
how to compactify the moduli space of Einstein metrics. 
The concept of a {\it Gromov-Hausdorff convergence}
provides a general framework for studying metric degenerations.  
Let $(M_k,p_k)$ is a sequence of complete pointed $n$-dimensional
Riemannian manifolds that Gromov-Hausdorff converge to the space $(Y,q)$;
in other words $(M_k,p_k)$ degenerates to $(Y,q)$. 
Due to {\it Gromov's compactness theorem}, a simple
way to ensure that $(M_k,p_k)$ has a Gromov-Hausdorff converging
subsequence is to assume that $\Ric(M_k)\ge c$ for some $c$. 
The sequence $M_k$ is said to {\it collapse near the
points $p_k$} if the volumes of the unit balls centered at $p_k$ 
tend to zero as $k\to\infty$. Otherwise, $M_k$ is called
{\it  non-collapsing near $p_k$}. 
An additional challenge is that collapsing and
non-collapsing may occur at the same time on different parts 
of the manifold (even though if $\Ric (M_k)\ge c$, then 
the distance between collapsed and non-collapsed parts 
of $M_k$ has to go to infinity as $k\to \infty$ because
of Bishop-Gromov's volume comparison).
It is mainstream of global Riemannian geometry
to study collapsing and non-collapsing sequence of manifolds
under various assumptions on curvature such as $|sec|\le C$, or
$\sec\ge 0$, or $\Ric\ge c$, or $|\Ric|\le c$, or 
various integral curvature bounds, or assuming that the metric is
Einstein, or K\"ahler-Einstein, or Calabi-Yau etc. 

\section{Manifolds with two sided bounds on sectional curvatures.}

To date the most complete picture of collapse 
is established for manifolds 
with two sided bounds on sectional curvature, i.e. 
when $|\seccurv(M_k)|\le c$. 
This bound is equivalent to the upper bound on the norm 
of the Riemann curvature tensor, hence it may seem
overly restrictive, especially when compared with bounds
on Ricci curvature which are much more natural in the 
context of {\it  general relativity}.
However, studying collapse under sectional curvature bounds,
rich and complex as it is, provides deep insights as to what
might happen under bounds on Ricci curvature.

If $|\seccurv(M_k)|\le c$, then 
any non-collapsing sequence $(M_k, p_k)$ has a 
subsequence that Gromov-Hausdorff converges to $(Y,q)$, where 
$Y$ a smooth manifold $Y$ with a $C^{1,\alpha}$-Riemannian metric, and
for any compact domain $D\subset Y$ there are
compact domains $D_k\subset M_k$, and diffeomorphisms
$\phi_k\co D\to D_k$ such that after pullbacking the metrics
from $D_k$ to $D$ via $\phi_k$,
the metrics converge in $C^{1,\alpha}$--topology  
(meaning that in each coordinate chart the coefficients of the
metric tensors $g_{ij}$ of $M_k$ converge in $C^{1,\alpha}$--topology).
If $D$ is Hausdorff close to the metric ball $B(q,R)$,
then $D_k$ is Hausdorff close to the metric ball $B(p_k,R)$.
If $M_k$ have uniformly bounded diameters, then one can take
$D_k=M_k$, in particular, for each positive $c,d, v$,
there only finitely many diffeomorphism
classes of complete Riemannian manifolds with $|\seccurv(M_k)|\le c$,
$\diam(M_k)\le d$, $vol(M_k)\ge v$.
(This result is sometimes called 
{\it Cheeger-Gromov's compactness theorem}.
Following a somewhat different route
Anderson~\cite{And-precomp} 
obtained the same conclusion under the weaker 
assumptions on Ricci curvature and
injectivity radius: $|\Ric(M_k)|\le c$ and $\inj_{p_k}(M_k)\ge i>0$).

The collapsing case with $|\seccurv(M_k)|\le c$ is much more complex.
The cornerstone of the theory is Gromov-Ruh's theorem on 
almost flat manifolds: a manifold $M$ admits a sequence of  
complete Riemannian metrics $g_k$ with $\diam (M,g_k)\le 1$ and
$|\seccurv (M,g_k)|\to 0$ as $k\to\infty$ if and only if $M$ is an
infranilmanifold~\cite{Gro-af, BusKar-af, Ruh-af}. 
An {\it infranilmanifold} is the quotient
of $N$ by a cocompact discrete torsion-free subgroup of
the semidirect product 
$N\rtimes F$, where $N$ connected simply-connected
nilpotent Lie group (such as $\mathbb R^m$, or the Heisenberg group),
and $F$ is a finite group of Lie group automorphisms of $N$.
For example, if $N=\mathbb R^m$, then infranilmanifolds
are precisely compact flat manifolds. Another important example is
a {\it nilmanifold}, i.e. the quotient of $N$
by its cocompact discrete subgroup. Topologically, nilmanifolds can be
characterized as iterated principle circle bundles, such as
a circle, a $2$-torus, all circle bundles over tori etc., and
any infranilmanifold is finitely covered by a nilmanifold.

With this definition 
almost flat manifolds are the solutions
of bounded size perturbations of the equation $\sec=0$, hence the name.
After suitable rescaling the assumptions
$\diam\le 1$, $|\seccurv|\to 0$ turn into $\diam\to 0$, 
$|\seccurv|\le 1$ so Gromov-Ruh's theorem characterizes 
manifolds that Gromov-Hausdorff converge (in fact collapse) 
to a point as infranilmanifolds. 
For example, any compact flat manifold collapses to a point simply 
by scaling the metric by $\epsilon$ with $\epsilon\to 0$.
A typical collapse 
of an infranilmanifold to a point necessarily involves inhomogeneous
scaling, e.g. for the circle bundle over the torus the fiber
could be scaled by $\epsilon^2$ while the base is scaled by $\epsilon$.
Thus the collapse generally takes place on several scales, and
an inhomogeneous scaling is essential to ensure that
the sectional curvature stays bounded during the collapse, 
i.e. that $|\seccurv|\le 1$.

{\it Fukaya's fibration theorem}~\cite{Fuk2006}
says that if $M_k$ collapses near $p_k$
to a Riemannian manifold $(Y,q)$, then for any compact domain $D\subset Y$
there are compact domains $D_k\subset M_k$
such that $D_k$ smoothly fibers over $D$ with infranilmanifolds
as fibers, and the collapse happens along these infranilmanifolds.
If $D$ is Hausdorff close to the metric ball $B(q,R)$,
then $D_k$ is Hausdorff close to the metric ball $B(p_k,R)$.
The case when $Y$ is not a manifold can be understood via
the equivariant version of Fukaya's fibration theorem which says
that the (orthonormal) frame bundles $FD_k$
with their natural metrics fiber over a certain Riemannian 
manifold $L$ with infranilmanifolds as fibers, the collapse 
occurs along the fibers, and the fibration is 
$O(n)$--equivariant with respect to some isometric 
$O(n)$--action on $L$ such that $L/O(n)$ is identified
with a subspace of $Y$. Pushing down
the fibration $FD_k\to L$ by the $O(n)$-action
one gets a stratification on $D_k$ whose strata 
are infranilmanifolds of various dimensions, and $D_k$
collapses to $L/O(n)\subset Y$ along the strata.

This description of local collapse is 
adequate for many applications, yet a much more general and
detailed picture was obtained by Cheeger-Fukaya-Gromov~\cite{CFG}. 
For any $\epsilon>0$, a complete Riemannian $n$-dimensional 
manifold with $|\seccurv|\le c$ can be partitioned into
two disjoint sets where the injectivity radius is $\ge\epsilon$
and $\le\epsilon$, called {\it $\epsilon$-thick} 
and {\it $\epsilon$-thin} part of the manifold. 
On compact domains in the $\epsilon$-thick
part one has Cheeger-Gromov compactness as explained above.
Results of Cheeger-Fukaya-Gromov describe the $\epsilon(n)$-thin
part for some  universal constant $\epsilon (n)$ depending 
only on $n$. They prove that 
the $\epsilon(n)$-thin part carries the so-called  
{\it invariant Riemannian
metric} (that can be chosen arbitrary close 
in $C^{1,\alpha}$-topology to the original metric) 
such that every point of the $\epsilon(n)$-thin
part has a neighborhood of definite size that is a tubular neighborhood of 
an embedded infranilmanifold, and there is an isometric action of a nilpotent 
Lie group in the universal cover of the neighborhood that stabilizers
the preimage of the infranilmanifold. Orbits of the
action descend to infranil strata in the $\epsilon(n)$-thin
part of the manifold. The family of local actions forms
what is called an {\it N--structure}, where ``N'' 
stands for ``nilpotent''. 

An N--structure captures all collapsed directions on every 
possible scale, yet it may be hard
to analyze partly because it carries so much information.
For many applications it is easier to deal with the so-called 
F--structures, where ``F'' stands for ``flat'', which were
introduced and studied by Cheeger-Gromov~\cite{CheGro-fI}. 
An F-structure capture information on collapse at
the smallest possible scale, that of the injectivity radius. 
In retrospect,
an F-structure corresponds to the center of an N--structure. 
An F--structure can be described as 
a collection of compatible local actions of tori
in finite covers of coordinate charts, so this notion generalizes
that of a torus action. By contrast, for the N--structures
the nilpotent group actions generally exist only
in the universal covers 
(which are almost never finitely-sheeted) of
``coordinate charts'' of the N-structure, and this 
makes N-structures harder to deal with.

The first nontrivial example of collapsing was studied by Berger
who considered the Hopf fibration $S^1\to S^3\to S^2$ where $S^3$
carries a constant curvature metric. Berger noted that multiplying
the metric in the fiber direction by $\epsilon$ while leaving
the metric in the orthogonal direction unchanged keeps the sectional
curvature bounded as $\epsilon\to 0$, so eventually $S^3$
looks like $S^2$ with a metric of $\seccurv=4$. By O'Neill's formula
this phenomenon extends to any principal torus bundles, and in fact to
any manifold that admits a torus action with no fixed point. Much more 
generally, Cheeger-Gromov~\cite{CheGro-fII} 
proved that any manifold $M$ carrying an 
F--structure of positive rank admits a sequence of metrics 
for which $(M,p)$ collapses with $|\seccurv|\le c$, 
where $p$ is an arbitrary point of $M$, and
``positive rank'' means that 
the local actions of tori have no fixed points.
Examples of compact manifolds that carry no 
F--structure of positive rank are those
of non-zero Euler characteristic (Cheeger-Gromov~\cite{CheGro-fII})
and those of nonzero simplicial volume (Paternain-Petean~\cite{PatPet}).
Thus manifolds satisfying these topological assumptions
cannot collapse with $|\seccurv|\le c$.

\section{Manifolds with lower bounds on sectional curvatures.}

By Gromov's compactness theorem the class of $n$-dimensional
complete Riemannian manifolds of $\seccurv\ge c$
is precompact in pointed Gromov-Hausdorff topology. The closure of the
class consists of the so-called Alexandrov spaces of dimension $\le n$ 
and curvature $\ge c$ in the comparison sense. A great deal
is known about geometry and topology of Alexandrov 
spaces~\cite{BBI, perelman-icm, Fuk2006}.
For the purposes of this article 
an {\it Alexandrov space} is a complete finite-dimensional
path metric space of curvature bounded below (where
a lower curvature bound is understood in
comparison sense, i.e. $\curv\ge c$ means that small triangles
in the space are at least as thick as the triangles in the 
plane of constant curvature $c$). 
Most notions of dimension
are equivalent for Alexandrov spaces, in particular, topological
and Hausdorff dimensions coincide. Any Alexandrov space contains 
an open dense set of points which have neighborhoods bi-Lipschitz equivalent
to Euclidean balls. Thus some analysis can be done
on Alexandrov spaces.
Alexandrov spaces admit stratifications into topological manifolds,
and every point has a (contractible) neighborhood
homeomorphic to the {\it tangent cone} at the point.

Every $2$-dimensional Alexandrov space is a manifold 
(possibly with boundary). Any convex subset of a Riemannian manifold of 
$\seccurv\ge c$ is an Alexandrov space of $\curv\ge c$.
Another good example of an Alexandrov space of $\curv\ge c$ is the 
quotient of a complete Riemannian manifold $M$ of $\seccurv\ge c$
by a compact isometry group $G$ of $M$. After giving $G$
a biinvariant metric of $\seccurv\ge 0$, 
the diagonal $G$-action on $\epsilon G\times M$ becomes
isometric, so $(\epsilon G\times M)/G$ carries the Riemannian 
submersion metric $g_\epsilon$
of $\seccurv\ge\min\{c,0\}$. The manifold
$(\epsilon G\times M)/G$ is diffeomorphic 
to $M$ with $g_\epsilon$-diameters of $G$-orbits converging to zero
as $\epsilon\to 0$. Therefore, $(M,g_\epsilon )$ 
collapse to $M/G$ as $\epsilon\to 0$
while $\seccurv (M,g_\epsilon ) \ge\min\{c,0\}$.

{\it Perelman's stability theorem}~\cite{Per-stab} 
says that if a non-collapsing sequence
of pointed complete Riemannian $n$-manifolds $(M_k,p_k)$ Gromov-Hausdorff
converges to the Alexandrov space
$(Y,q)$, then for any compact domain $D\subset Y$ and for
all large $k$ there are compact domains $D_k\subset M_k$, and
homeomorphisms $h\co D\to D_k$ satisfying 
$|d(x,y)-d_k(h(x),h(y))|\to 0$
for all $x,y\in D$ as $k\to\infty$. 
If $D$ is Hausdorff close to the metric ball $B(q,R)$,
then $D_k$ is Hausdorff close to the metric ball $B(p_k,R)$.
In particular, 
if $\diam(M_k)\le d$,
then $M_k$ is homeomorphic to $Y$ for all large $k$.

{\it Yamaguchi's fibration theorem}~\cite{Yam-fib} says that
if $(M_k,p_k)$ is a collapsing sequence 
pointed complete Riemannian $n$-manifolds with 
$\seccurv\ge c$ that Gromov-Hausdorff 
converges to a Riemannian manifold $(Y,q)$, then 
for any compact domain $D\subset Y$
there are compact domains $D_k\subset M_k$
such that $D_k$ smoothly fibers over $D$ such that 
the fibers have almost nonnegative curvature (in a certain weak sense).
If $D$ is Hausdorff close to the metric ball $B(q,R)$,
then $D_k$ is Hausdorff close to the metric ball $B(p_k,R)$.
Burago-Gromov-Perelman~\cite{BGP} and Yamaguchi~\cite{Yam-fib-alex} 
also proved a version of this result when $Y$ 
is an Alexandrov space with metrically mild singularities.
This fibration theorem may seem very special, yet it allows
for inductive reasoning, and it can be used repeatedly
in combination with other methods such as rescalings and splitting, 
to yield e.g. the following result of Fukaya-Yamaguchi~\cite{FukYam-p1}: 
if $M_k$ Gromov-Hausdorff
converges to a point with $\seccurv\ge c$, then
for large $k$ the fundamental group of $M_k$ has a nilpotent
subgroup of finite index.

\section{Manifolds with various bounds on Ricci curvature.}

By Gromov's compactness theorem any sequence of pointed
complete Riemannian $n$-manifolds  satisfying $\Ric\ge c$ 
has a subsequence $(M_k,p_k)$ that Gromov-Hausdorff converges to a 
locally compact complete metric space $(Y,q)$
of Hausdorff dimension $\le n$. It is generally thought that away
from a small singular set the space $Y$ should
look like a Riemannian manifold and the corresponding parts of 
$M_k$ converge/collapse to the complement of the singular set 
in the manner similar to what happens under sectional
curvature bounds. 

In the non-collapsing case this kind of results
were obtained by Cheeger-Colding as described below 
(see~\cite{Chi-Ric-book}). 
A point of $Y$ is called {\it regular} if each of its tangent 
cones is isometric to the Euclidean space of some dimension; 
denote by $\mathcal R$ 
the set of regular points. Then $\mathcal S:=Y\setminus R$ 
denotes the set of {\it singular} points. 
Let $\mathcal R_\epsilon$ be the set of points of $Y$ such
that for every $y\in \mathcal R_\epsilon$ each  
tangent cone at $y$ has the property that the Gromov-Hausdorff
distance between the unit ball in $\mathbb R^n$
and the unit ball centered at the apex of the tangent
cone is $<\epsilon$. Clearly, 
$\mathcal R\subset\mathcal R_\epsilon$. 

Assuming that $\Ric(M_k)\ge c$, and that the sequence 
$(M_k,p_k)$ is non-collapsing
and Gromov-Hausdorff converges to $(Y,q)$, Cheeger-Colding
proved the following. 
The space $Y$ has Hausdorff dimension $n$,
and the isometry group of $Y$ is a Lie group.
The set of regular points $\mathcal R$ is connected, dense, and 
each tangent cone of $y\in \mathcal R$
is isometric to $\mathbb R^n$. The Hausdorff dimension
of $S$ is $\le n-2$, which is optimal as seen from 
$2$-dimensional examples. There exists $\epsilon (n)$
depending only on $n$ such that if $\epsilon\le \epsilon (n)$, 
then $\mathcal R$ lies in the interior of $\mathcal R_\epsilon$,
and furthermore, $\mathcal R_\epsilon$ 
is a connected smooth manifold on which
the metric is bi-H\"older to a smooth Riemannian metric. 

If the assumption $\Ric\ge c$ is replaced by
$|\Ric|\le c$, Cheeger-Golding 
proved that $R=R_\epsilon$ for 
$\epsilon\le \epsilon (n)$ which implies that $R$ is open 
and hence $S$ is closed, and
furthermore, $R$ is a $C^{1,\alpha}$-Riemannian manifold,
and convergence $M_k$ to $Y$ is $C^{1,\alpha}$-topology.
The same is true in $C^\infty$-topology if
each $M_k$ is Einstein, in which case $R$ is also Einstein.

Cheeger-Colding also obtained various 
results in the collapsing case 
(under the assumption $\Ric (M_k)\ge c$), e.g. the Hausdorff dimension
of $Y$ in the collapsing case is $\le n-1$, yet at present 
the picture is much less complete.
One of Cheeger-Colding's goals in the collapsing case was
to prove Gromov's conjecture that if $\Ric(M_k)\ge c$ and
$Y$ is a point, then $\pi_1(M_k)$ contains a 
nilpotent subgroup of finite index for all large $k$.
A major obstacle was
the absence of the fibration theorem 
in the Ricci curvature setting (indeed, Anderson
gave examples of almost Ricci-flat $4$-manifolds
that collapse to the $3$-torus but cannot fiber over it).
Recently, a solution of Gromov's conjecture was announced
by Kapovitch-Wilking.

A variety of results on degeneration of metrics with
Ricci curvature bounds can be obtained under additional 
bounds on the $L^p$-norm of the Riemann curvature tensor 
$||R||_{p}$. To illustrate, consider the convergence theorem
of Anderson for the class of $n$-dimensional compact 
Riemannian manifolds satisfying  the bounds
$|\Ric|\le c$, $\diam\le d$, $\Vol\ge v>0$, 
$||R||_{n/2}\le C$~\cite{And-surv}. 
Any sequence of manifolds in the class has a subsequence
converging to a Riemannian orbifold with finitely many 
singular points. The metrics converge in $C^{1,\alpha}$-topology
away from the singular set, and again if each manifold in the
sequence is Einstein, the convergence is in $C^{\infty}$-topology.
For example, a sequence of Eguchi--Hanson Ricci-flat 
metrics on the tangent bundle 
to the $2$-sphere converges in this way 
to a flat metric cone on $\mathbb{RP}^3$.

By Gauss-Bonnet's theorem the bound $||R||_{n/2}\le C$ comes for free
in studying Einstein metrics on a {\it fixed}
compact $4$-dimensional manifold. In this case the metric degeneration
was extensively studied by Anderson~\cite{And-surv}, and more recently, by
Cheeger-Tian~\cite{CheTia}, and roughly speaking, they prove that
away from a finite collection
of ``blowup points'' where the curvature concentrates,
the convergence/collapsing occurs similarly to what happens under
the two-sided bound on sectional curvature. 

In higher dimensions little is known on degenerations of Einstein 
manifolds. In fact, even the case of {\it Calabi-Yau manifolds} 
(with their Ricci-flat K\"ahler metric) is wide-open. It is worth 
mentioning that
degenerations of Calabi-Yau manifolds to
the so called large complex structure limit point
is of interest for the {\it mirror symmetry}.
If $M$ is a simply-connected Calabi-Yau $n$-fold,
a {\it large complex structure limit point} is a point in
the compactified moduli space of the complex structures on $M$
which in a sense represents the ``worst possible degeneration''
of the complex structure. There is a conjecture of 
Kontsevich-Soibelman and Gross-Wilson, 
that if $(M, g_k)$ are Ricci-flat K\"ahler 
metrics whose complex structure converge to a
large complex structure limit point, and if $\diam(M,g_k)$
is bounded away from zero and infinity as $k\to\infty$, 
then $(M, g_k)$ collapses along a singular torus fibration
to a metric space homeomorphic to the $n$-sphere 
(i.e. the dimension of the sphere is half of 
the real dimension of $M$). This is relevant for the mirror symmetry
because, conjecturally, the mirror manifold 
is obtained by dualising the fibration. 
Gross-Wilson~\cite{GrossWil} 
verified the conjecture for the $K3$-surface.

\appendix
\section{Gromov--Hausdorff distance}

Gromov--Hausdorff distance measures how far the abstract
metric spaces are from being isometric. This notion was
introduced by Gromov in~\cite{Gro} 
as a generalization of the classical 
notion of Hausdorff distance.
If $X$ and $Y$ are two non-empty compact
subsets of a metric space $S$, 
then the {\it Hausdorff distance} $d_H(X,Y)$ is the infimal 
number $r$ such that the closed $r$--neighborhood of $X$ 
contains $Y$ and the closed $r$--neighborhood of $Y$ contains $X$. 
If $X$ and $Y$ are 
two compact metric spaces, then the 
{\it Gromov-Hausdorff distance}
$d_{GH}$ between $X$ and $Y$ 
is defined to be the infimum of the numbers 
$d_H(f (X ), g (Y ))$ over all metric spaces $S$ and 
all distance-preserving 
embeddings $f\co X\to S$, $g\co Y\to S$.
Like any distance, the Gromov-Hausdorff distance defines 
a notion of convergence for sequences 
of compact metric spaces, called the
{\it Gromov-Hausdorff convergence}.
Up to isometry the limits are unique, because
$d_{GH}(X,Y)=0$ if and only if $X$ and $Y$ are isometric. 

A {\it pointed Gromov-Hausdorff convergence} is 
a version of Gromov-Hausdorff convergence suitable to deal with 
non-compact metric spaces. This is
analogous to uniform convergence of functions on compact subsets,
and a basepoint can be thought of as the position of the observer.
A sequence $(X_k, x_k)$ of 
pointed locally compact complete path metric spaces is said to 
{\it Gromov-Hausdorff converge to $(Y,y)$}  
if for any $R > 0$ the sequence
\[
\inf_{f,g,S}\{ d_H(f (B(x_k, R)), g (B(y,R )))+ d_H(f(x_k),g(y))\}
\]
converges to zero as $k\to \infty$,
where  $B(x_k, R)$, $B(y,R)$ are the closed $R$-balls
centered at $x_k$, $y$, and the infimum is taken
over all metric spaces $S$ and 
all distance-preserving embeddings $f\co X\to S$, $g\co Y\to S$.
In other words, for each $R$, the balls $B(x_k, R)$ 
Gromov-Hausdorff converge to $B(y,R)$ and furthermore
their centers converge as well.
The limit $Y$ is automatically a complete locally compact
path metric space (where ``path'' means that the distance
between any two points is the infimum of lengths of paths
joining the points). Again, the limits are unique up to
a basepoint-preserving isometry.

A family of spaces is called {\it Gromov-Hausdorff precompact} 
if any sequence in the family has
a convergence subsequence. A useful criterion of
Gromov-Hausdorff precompactness says that $(X_k,x_k)$
is Gromov-Hausdorff precompact if and only if for every 
$R$, $\epsilon$ there exists an integer $N(R,\epsilon)$
such that for each $k$ the ball
$B(x_k,R)$ can be covered by $N(R,\epsilon)$-balls
of radius $\epsilon$. 

Let $r X$ denote a copy of $X$
with the metric multiplied by a constant $r>0$. 
Given a sequence of numbers $\{r_k\}$ satisfying
$r_k\to\infty$ as $k\to \infty$,
the Gromov-Hausdorff limit of $(r_k X, p)$ 
(if it exists) is called the {\it tangent cone} of $X$ at $p$
defined by the sequence of scaling factors $\{r_k\}$.
This is an invariant of the local geometry of $X$ at $p$.
As a trivial example, the tangent cone of 
an $n$-dimensional Riemannian manifold at any point
is isometric to $\mathbb R^n$ which manifests
a familiar fact that Riemannian manifolds are infinitesimally 
Euclidean.

In the opposite direction, asymptotic properties 
of the space are often captured by the limit (if it exists)
of the sequence
$(\frac{1}{r_k} X,p)$ where $r_k\to\infty$ as $k\to \infty$
which is called an {\it asymptotic cone} of a metric space $X$
at $p$ defined by $\{r_k\}$.
Gromov used asymptotic cones to prove 
that every group of polynomial growth
has a nilpotent subgroup of finite index by taking $X$
to be a Cayley graph of the group equipped with the word metric.

In general, the Gromov-Hausdorff limit of $(X_k,x_k)$
need not exist, e.g. the {\it hyperbolic plane} has no
asymptotic cones (which is due to the fact that the volume of
balls in the hyperbolic plane grows exponentially rather than 
polynomially). Fortunately,
there is a generalization of the Gromov-Hausdorff
convergence (due to van den Dries-Wilkie)
that associates a limit to {\it any} sequence of pointed 
metric spaces. 
For the purposes of this article
a {\it non-principal ultrafilter} $\omega$
is a device that assigns a unique
limit $\lim_\omega a_k$
to any bounded sequence $\{a_k\}$ of real numbers.
No explicit examples non-principal ultrafilters are known,
yet they exist by Zorn's Lemma.
Now if $(X_k,p_k)$ is a sequence of pointed metric spaces,
one lets $X_\infty$ be the set of sequences $\{x_k\}$ with
$x_k\in X_k$
such that the distances $d_k(x_k,p_k)$ are uniformly bounded,
and one defines 
$d_\infty(\{x_k\},\{y_k\}):=\lim_\omega d_k(x_k,y_k)$. 
The space $(X_\infty,d_\infty)$
is a pseudometric space, and usually there are distinct 
points $\{x_k\},\{y_k\}$ with $d_\infty(\{x_k\},\{y_k\})=0$;
in this case one deems $\{x_k\},\{y_k\}$ equivalent, and
the set of equivalence classes becomes a metric space
which is called the {\it ultralimit} $X_\omega$ of $(X_k,p_k)$.
It is known that 
if a sequence $(X_k,p_k)$ Gromov-Hausdorff converges to $(Y,q)$, 
then the ultralimit of $(X_k,p_k)$ is isometric to $Y$. 

The ultralimit of the sequence 
$(\frac{1}{r_k} X,p)$, where $r_k\to\infty$
as $k\to \infty$, is also called an {\it asymptotic cone}
of $X$ at $p$. With this new definition, the asymptotic cone of a 
hyperbolic plane is a tree with 
uncountable branching at every point. Asymptotic cones of groups
have been especially useful in geometric group theory, 
and in studying quasi-isometries of groups and nonpositively 
curved spaces.

The notion of Gromov-Hausdorff convergence provides a
convenient framework for studying degenerations of Riemannian 
manifolds. A basic fact is the
{\it Gromov's compactness theorem}, 
which implies (via Bishop-Gromov's 
volume comparison) that for each number $c$
the set of Riemannian $n$-dimensional
manifolds with Ricci curvature $\Ric\ge c$  
is Gromov-Hausdorff precompact, and
any limit metric space has Hausdorff dimension $\le n$. 
The same result holds for unpointed Gromov-Hausdorff distance
when restricted to metric spaces of uniformly bounded diameter.

For more information on Gromov-Hausdorff distance
see~\cite{Gro, BH, BBI}.

\section{Acknowledgements}

Thanks are due to Vitali Kapovitch for helpful comments on 
the first version of this paper.
The author was partially supported by the NSF grant \# DMS-0352576. 

\providecommand{\bysame}{\leavevmode\hbox to3em{\hrulefill}\thinspace}
\providecommand{\MR}{\relax\ifhmode\unskip\space\fi MR }
\providecommand{\MRhref}[2]{%
  \href{http://www.ams.org/mathscinet-getitem?mr=#1}{#2}
}
\providecommand{\href}[2]{#2}

\end{document}